 \documentstyle[11pt]{article}  
   \author {Bruno Fabre \\
   22, rue Emile Dubois, 75014 Paris, France}
   
     \date{}
     \title {Sur la cohomologie de Dolbeault des vari\'et\'es projectives et les courants localement
     r\'esiduels}  
       
\begin{document}
          \def\dem{{\it Proof.}\\} 
          \def\CQFD{\hbox{\vrule height 5pt depth 5pt width 5pt}}   \def\bk{\backslash}  
             \def\B{ \mbox{I\hspace{-.15em}B}}
               \def\N{ \mbox{I\hspace{-.15em}N}} 
                \def\Z{ \mbox{Z\hspace{-.3em}Z}}  \def\R{\mbox{I\hspace{-.17em}R}}  
                \def\P{ \mbox{I\hspace{-.17em}P}} 
                 \def\C{ \mbox{l\hspace{-.47em}C}}
                   \def\Q{ \mbox{l\hspace{-.47em}Q}} 
                     \newtheorem {ethm}{Th\'eor\`eme \hspace{2pt}} 
                       \newtheorem {thm}{Theorem\hspace{2pt}} 
\newtheorem {ecor}{Corollary\hspace{2pt}} 
 \newtheorem {cor}{Corollaire \hspace{2pt}} 
\newtheorem {edefi}{Definition\hspace{2pt}} 
\newtheorem {defi}{D\'efinition\hspace{2pt}} 
  \newtheorem {elem}{Lemma\hspace{2pt}}
  \newtheorem {lem}{Lemme\hspace{2pt}}

\maketitle

{\small
{\bf R\'esum\'e.}

Soit $X$ une vari\'et\'e projective.
Soit $Y_1,\dots,Y_{p+1}$ des hypersurfaces sur $X$, en position d'intersection compl\`ete, et associ\'ees \`a des diviseurs amples.
Nous montrons dans cette note que pour $i\le p$, la cohomologie de Dolbeault $H^i(\Omega^q)$ du faisceau $\Omega^q$ des $q-$formes holomorphes sur $X$ peut \^etre calcul\'ee comme la cohomologie de degr\'e $i$ d'un complexe de sections globales de courants localement r\'esiduels. 
On peut \'egalement calculer de la m\^eme mani\`ere la cohomologie du sous-faisceau ${\tilde\Omega}^q$ des formes holomorphes $\partial-$ferm\'ees comme la cohomologie d'un sous-complexe de courants localement r\'esiduels $\partial-$ferm\'es.
On en d\'eduit comme corollaire que toute classe de type $(i,i)$ $(i\le p+1)$ admet comme repr\'esentant un courant $d-$ferm\'e \`a support dans $Y:=Y_1\cap \dots\cap Y_i$.
On montre \'egalement que tout courant localement r\'esiduel $T$ sur $X$, \`a support dans l'intersection $Y=Y_1\cap \dots\cap Y_{i-1}$, peut \^etre \'ecrit comme un r\'esidu global $T=Res_{Y_1,\dots,Y_{i-1}}{\Psi}$,
avec $\Psi$ une $q-$forme m\'eromorphe \`a p\^ole dans $Y_1\cup\dots\cup Y_i\cup Y_{i}$, et on peut se passer de $Y_{i}$ ssi $T$ est $\overline\partial-$exact.
On retrouve ainsi un th\'eor\`eme de Hererra-Dickenstein-Sessa ({\cite{HDS}). On donne pour conclure une nouvelle formulation \'equivalente de la conjecture de Hodge.

{\bf Abstract.} 
Let $X$ be a projective manifold. Let $Y_1,\dots,Y_{p+1}$ be $p+1$ ample hypersurfaces in complete intersection position on $X$, each defined by the global section of an ample Cartier divisor. We show in this note that for $i\le p+1$, the cohomology groups $H^i(\Omega^q)$ can be computed as the $i-$th cohomology groups of some complex of global sections of locally residual currents on $X$. We could also compute the cohomology of the subsheaves ${\tilde\Omega}^q\subset \Omega^q$ of $\partial-$closed holomorphic forms by the corresponding subsheaves of $\partial-$closed locally residual currents. We deduce like this that any cohomology class of bidegree $(i,i)$ has an element which is a $d-$closed  locally residual current with support in $Y_1\cap \dots\cap Y_i$.
We also show that any locally residual current $T$ of bidegree $(q,i-1)$ with support in $Y_1\cap \dots Y_{i-1}$ can be written as a global residue $T=Res_{Y_1,\dots,Y_{i-1}}{\Psi}$ of some meromorphic form with pole in $Y_1\cup\dots\cup Y_{i}$. We can avoid $Y_{i}$ iff the current in $\overline\partial-$exact; we deduce as corollaries a theorem of Hererra-Dickenstein-Sessa ({\cite{HDS}). 
We give as a conclusion a new formulation of the Hodge conjecture.
}

\section{Abridged English version}

Let ${\cal G}$ be a sheaf on a topological space $X$. We know from \cite{Har} the existence
for any two closed subsets $Y,Y' \subset X$, of the functor $\Gamma_{Y\backslash Y'}$, and the right derived functors ${\cal H}^i_{Y\backslash Y'}$; and moreover,
 for any three closed subsets $Y,Y',Y''\subset X$, of a {\it topological residue operator}:

$$Res_Y: {\cal H}^i_{Y'\backslash {Y\cup Y''}}({\cal G})\to {\cal H}^{i+1}_{Y'\cap Y\backslash Y''}({\cal G})$$ 

Let now $X$ be a complex analytic manifold of pure dimension $p$.
We have as a natural topology on $X$ the Zariski topology.
Let us consider the sheaves for the the Zariski topology, as $\C_X-$modules, where $\C_X$ is the constant sheaf. 

We consider a subcategory of sheaves ${\cal F}$ satisfying the two following properties:
i) For any Zariski-closed subset $Y$, if there is a section $s\in {\cal F}_x$, then there is an element $f\in {\cal O}_x$ and an integer $k$ such that $f^k s=0$ in ${cal F}_x$;

ii) If $s$ is a section in some Zariski open subset $U\subset X$, then for any $x\in X\backslash U$, there is an element $f\in {\cal O}_x$ and an integer $k$ such that $f^k s$ extends to an element of ${\cal F}_x$. 

In this category of Zariski sheaves, any sheaf ${\cal F}$ is contained in a flabby sheaf, namely his associated meromorphic sheaf $\tilde{\cal F}$, defined as follows: a section $s\in \tilde{\cal F}(U)$ is a section $s'\in {\cal F}(V)$, for $V$ some dense Zariski open subset of $U$. 
So we can define the derived functors of $\Gamma_Y$ (with $Y$ locally closed) in this category, and this associate to a Zariski sheaf $\cal F$ in this category a sheaf ${\cal H}_Y({\cal F})$.

For the sheaves ${\cal G}$ in the initial category of ${\cal O}_X-$modules, we associate a Zariski sheaf ${\cal G}'$, defined for a Zariski open subset $U$ by defining ${\cal G}'(U)$ as the sections of ${\cal G}(U)$ which extends meromorphically outside $U$. Then, if ${\cal G}'$ belongs to the above subcategory of Zariski-sheaves, we define the {\it moderate cohomology} as 
${\cal H}^i_{[Y]}({\cal G}):={\cal H}^i_Y({\cal G}')$.

From this, we deduce for these sheaves $\cal G$ a moderate residue operator:
$$Res_Y: {\cal H}^i_{[Y'\backslash {Y\cup Y'']}}(U,{\cal G})\to {\cal H}^{i+1}_{[Y'\cap Y\backslash Y'']}(U,{\cal G}).$$

Let us now assume that ${\cal G}=\Omega^q$.

We have thus a long exact sequence of moderate cohomology groups:
$$\begin{array}{c} 
0\to H^0_{[Z\cap Y]}(\Omega^q)\to H^0_Z(\Omega^q)\to H^0_{[Z\backslash Y]}(\Omega^q) \\ \to 
H^1_{[Z\cap Y]}(\Omega^q)\to 
\dots \to H^i_{[Z\backslash Y]}(\Omega^q)\to H^{i+1}_{[Z\cap Y]}(\Omega^q)\to \dots
\end{array}$$

Let $Y$ an analytic hypersurface, and ${\cal F}$ any ${\cal O}_X-$module on $X$; we denote ${\cal F}(\star Y)$ the subsheaf of $\Gamma_U({\cal F})$, whose sections extend to a section of ${\cal F}$ if we multiply by some regular function vanishing on $Y$. 

\begin{ethm}[S. Boucksom]
Let us assume that $Z$ is a closed analytic subset of pure codimension $p$, and $Y$ an analytic hypersurface intersecting $Z$ properly.
Then, the moderate cohomology group $H^p_{[Z\backslash Y]}(\Omega^q)$ is isomorphic to the group ${\cal C}_Z(\star Y)$ of locally residual currents with support in $Z$, $\overline\partial-$closed outside $Y$; moreover,
the residue operator ${Res_Y}': {\cal C}^{q,i}_Z(\star Y)\to {\cal C}^{q,i+1}_{Z\cap Y}$
coming from the above $Res_Y: {\cal H}^i_{[Z\backslash Y]}(\Omega^q)\to {\cal H}^{i+1}_{[Z\cap Y]}(\Omega^q)$
is precisely the $\overline\partial-$operator.
\end{ethm}

Now, let $Y_1,\dots,Y_p$ analytic hypersurfaces of $X$ in complete intersection position.
We have a composed residue operator:
$Res_{Y_1,\dots,Y_k}: H^i_{Y\backslash Y'}({\cal G})\to H^{i+k}_{Y_1\cap \dots\cap Y_k\cap Y\backslash Y''}({\cal G}).$

Then, we have the following exact sequence of sheaves:

$$0\to \Omega^q\to {\cal C}^{q,0}(\star Y_1)\to \dots \to {\cal C}^{q,p}_{Y_1\cap\dots\cap Y_p}\to 0,$$
where the operator are $\overline\partial$.

Let us call {\it positive} an hypersurface $Y$ such that, for any coherent sheaf $\cal F$ on $X$, we have for $j>0$ $H^j({\cal F}(\star Y))=0$.

We have the following:

\begin{ethm}
For $i\le p+1$, the sheaves ${\cal C}^{q,i-1}_{Y_1\cap\dots\cap Y_{i-1}}(Y_{i})$ are acyclic, so that we have a canonical isomorphism:
$$H^i(\Omega^q)\simeq H^0(   {\cal C}^{q,i}_{Y_1\cap \dots\cap Y_i})/\overline\partial(H^0( {\cal C}^{q,i-1}_{Y_1\cap\dots\cap Y_{i-1}}))$$
\end{ethm}

Moreover, the elements of $H^0( {\cal C}^{q,i-1}_{Y_1\cap\dots\cap Y_{i-1}})$ can be represented as global residues, by the following theorem:
\begin{ethm}
If $T$ is a locally residual current with support in $Y_1\cap \dots\cap Y_{i-1}$, $\overline\partial-$closed outside $Y_{i}$, then $T$ can be written as a global residue:
$T=Res_{Y_1,\dots,Y_{i-1}}{\Psi}$, with $\Psi$ a meromorphic form with poles contained in $Y_1\cup\dots\cup Y_{i}$. Moreover, $T$ is $\overline\partial-$exact iff we can assume that $\Psi$ has his poles contained in $Y_1\cup\dots\cup Y_{i-1}$.
\end{ethm}

The proofs contains as principal ingredient that the sheaves ${\cal C}^{q,i}_{Y_1\cap\dots\cap Y_i}(\star Y_{i+1})$ are acyclic.

The assertion in the last theorem, saying that if the locally residual current is $\overline\partial-$exact, it can be written as a global residue, was already a theorem of Hererra-Dickenstein-Sessa (\cite{HDS}).

Let us define the following subcomplex:
$$0\to {\tilde \Omega}^q\to {\tilde{\cal C}}^{q,0}(\star Y_1)\to \dots \to {\tilde{\cal C}}^{q,p}_{Y_1\cap \dots \cap Y_p}\to 0.$$
where the subsheaves ${\tilde{\cal C}}^{q,i}_{Y_1\cap \dots\cap Y_i}\subset {{\cal C}}^{q,i}_{Y_1\cap\dots\cap Y_i}$ are defined by the $\partial-$closed condition. Then the same reasonning shows that this subcomplex computes the cohomology groups $H^i({\tilde\Omega}^q)$.

When $X$ is projective, the positivity condition on the $Y_i$ are satisfied, as soon are the associated line bundles are ample, by Serre's vanishing theorem.
Thus, we get that any $\overline\partial-$cohomology class of bidegree $(q,i)$, admits as representative a $d-$closed locally residual current with support in $Y_1\cap \dots\cap Y_i$. 

Let us recall that Hodge's conjecture is equivalent to saying that a cohomology class of bidegree $(i,i)$, which is moreover {\it integral}, which means that his integral over any real subvariety of dimension $2i$ is an integer, and such that his integral over all the complex subvarieties of dimension $i$ is zero, is zero. By the above theorem, it is equivalent to saying that an {\it integral} $d-$closed residual current with support on $Y_1\cap \dots\cap Y_i$, which is globally residual on any complex $i-$dimensional subvariety, is globally residual on $X$.

{{\bf Remark.}
For $q=n$, can can restrict to simple poles in the complex 
$$0\to \Omega^n\to {\cal C}^{n,0}(\star Y_1)\to \dots \to {\cal C}^{n,p}_{Y_1\cap\dots\cap Y_p}\to 0.$$
Then, by Kodaira's vanishing theorem, we get that this subcomplex also computes Dolbeault cohomology groups $H^i(\Omega^n)$. This is a variant of the main theorem of \cite{Rosly}.
}

\section{Enonc\'e des r\'esultats}

Soit $X$ un espace topologique; on connait d'apr\`es \cite{Harthorne}, pour deux ferm\'es $Y,Y'$, le foncteur $\Gamma_{Y\backslash Y'}$ ainsi que les foncteurs d\'eriv\'es ${\cal H}^i_{Y\backslash Y'}$.

On a une suite exacte:
$$\Gamma_Y {\cal F}\to {\cal F}\to \Gamma_{X\backslash Y}{\cal F},$$
avec un z\'ero \`a droite lorsque ${\cal F}$ est flasque. De plus, si ${\cal F}$ est flasque, il en est de m\^eme de $\Gamma_{Y'\backslash Y''}{\cal F}$.

Consid\'erons donc une r\'esolution flasque de ${\cal F}$:
$$0\to{\cal F}\to {\cal I}^0\to {\cal I}^1\to \dots $$
On en d\'eduit une suite exacte courte de complexes:
$$0\to \Gamma_{Y\cap (Y'\backslash Y'')}{\cal F}\to \Gamma_{Y'\backslash Y''}{\cal F}\to \Gamma_{Y'\backslash (Y\cup Y'')}{\cal F}\to 0$$,
d'o\`u on  d\'eduit par le lemme du serpent une suite exacte longue de cohomologie; on appelle les op\'erateurs:
$$Res_Y: {\cal H}^i_{Y'\backslash {Y\cup Y''}}({\cal F})\to {\cal H}^{i+1}_{Y'\cap Y\backslash Y''}({\cal F})$$ 
intervenant dans cette suite les op\'erateurs de {\it r\'esidus cohomologiques}.

Soit maintenant $X$ une vari\'et\'e analytique de dimension $n$.
On a sur $X$ une topologie naturelle de Zariski.
On a dans la cat\'egorie des faisceaux pour cette topologie de Zariski une sous-cat\'egorie, d\'efinie par les deux conditions suivantes:
1. Une section \`a support dans un sous-ensemble analytique, doit \^etre annul\'ee par une puissance de l'id\'eal ${\cal I}_Y$;
2. Une section sur un ouvert de Zariski $U=X\backslash Y$, s'\'etend \`a travers $Y$ apr\`es une multiplication par un \'el\'ement de ${{\cal I}_Y}_x$.
Dans cette cat\'egorie, tout faisceau ${\cal F}$ est inclus dans un faisceau flasque, en effet le faisceau $\tilde{\cal F}$ obtenu en prenant les sections $\tilde{\cal F}(U)$ comme les sections ${\cal F}(V)$ pour un ouvert de Zariski $V\subset U$.

On associe \`a un ${\cal O}_X-$module ${\cal G}$ un faisceau de Zariski ${\cal G}'$, en prenant sur un ouvert de Zariski $U$ les sections de ${\cal G}(U)$ qui se prolongent m\'eromorphiquement en dehors de $U$. On suppose que $X$ admet un recouvrement fini par des ouverts de Zariski de Stein.
Pour un ${\cal O}_X-$module ${\cal G}$ tel que ${\cal G}'$ appartient \`a la sous-cat\'egorie des faisceaux de Zariski d\'efinie ci-dessus, on d\'efinit la cohomologie mod\'er\'ee \`a support comme:
$${\cal H}^i_{[Y\backslash Y']}({\cal G}):= H^i_{Y\backslash Y'}({\cal G}').$$

On en d\'eduit comme ci-dessus, pour un tel ${\cal O}_X-$module, une suite exacte longue de cohomologie mod\'er\'ee, avec des op\'erateurs de r\'esidus:
$$Res_Y: {\cal H}^i_{[Y'\backslash {Y\cup Y''} ]}({\cal G})\to {\cal H}^{i+1}_{[Y'\cap Y\backslash Y'']}({\cal G})$$ 

Pour un sous-ensemble analytique ferm\'e $Y$, on note ${\cal F}(\star Y)$ les sections de ${\cal F}$ qui s'\'etendent m\'eromorphiquement \`a travers $Y$.

D'autre part, un courant localement r\'esiduel de bidegr\'e $(q,p)$ est un courant qui s'exprime localement sous la forme:
$\omega\wedge [1/f_{p+1}] \overline\partial [1/f_1]\wedge \dots \wedge [1/f_p]$,
avec $\omega$ une $q-$forme holomorphe et $(f_1,\dots,f_{p+1})$ une suite r\'eguli\`ere de fonctions holomorphes.

On a le lemme suivant:
\begin{lem}[S. Boucksom]
Soit $Y$ un sous-ensemble analytique de codimension pure $p$, et $Y'$ une hypersurface de $X$ coupant $Y$ proprement; les courants localement r\'esiduels de bidegr\'e $(q,p)$, \`a support dans $Y$ et $\overline\partial-$ferm\'es
(resp. $\overline\partial-$ferm\'es en dehors de $Y'$), forment un faisceau ab\'elien,
not\'e ${\cal C}^{q,p}_Y$ (resp. ${\cal C}_Y(\star Y') $).
\end{lem}

Soit $Y$ un sous-ensemble analytique de codimension pure $p$, et $Y'$ une hypersurface intersectant $Y$ proprement.

 Alors:

\begin{thm}\label{th0}
On a un isomorphisme canonique: 
$$H^p_{[Y\backslash Y']}(\Omega^q)\simeq {\cal C}^{q,p}_Y(\star Y')$$
De plus, \`a travers cet isomorphisme, l'op\'erateur ci-dessus
$$Res_{Y'}: H^p_{[Y\backslash Y']}(\Omega^q)\to  H^{p+1}_{[Y\cup Y']}(\Omega^q)$$
s'exprime comme $\overline\partial$ des courants.
\end{thm}

Plus g\'en\'eralement, si $Y''$ est une autre hypersurface coupant $Y\cap Y'$ proprement, on obtient un op\'erateur de r\'esidu $Res_{Y'}'$ sur les courants localement r\'esiduels \`a support dans $Y$ et $\overline\partial-$ferm\'es en dehors de $Y'\cup Y''$, en prenant l'op\'erateur correspondant au morphisme de r\'esidu:
$$Res_{Y'}: {\cal H}^p_{[Y\backslash (Y'\cup Y'')]}(\Omega^q)\to  {\cal H}^{p+1}_{[Y\cup Y'\backslash Y'']}(\Omega^q)$$
\`a travers l'isomorphisme ci-dessus.
Soit donc $Y_1,\dots, Y_p$ $p$ hypersurfaces en position d'intersection compl\`ete; on obtient les r\'esidus compos\'es par la compositon des op\'erateurs r\'esidus:
$$Res_{Y_1,\dots,Y_p}'=Res_{Y_1}'\circ \dots\circ Res_{Y_p}'$$

Alors, on en d\'eduit:
\begin{lem} Pour $\Psi$ une forme m\'eromorphe \`a p\^ole dans $Y_1\cup\dots\cup Y_p$, on a:
$Res_{Y_1,\dots,Y_p}'=\overline\partial Res_{Y_1,\dots,Y_{p-1}}'$.
\end{lem}

De plus, on a:

\begin{lem}
Le complexe suivant:
$$0\to \Omega^q\to {\cal C}^{q,0}(\star Y_1)\to \dots \to {\cal C}^{q,p}_{Y_1\cap\dots\cap Y_p}\to 0,$$
avec $\overline\partial$ comme op\'erateurs, est exact.
\end{lem}

En effet, un courant $\overline\partial-$ferm\'e s'\'ecrit localement comme r\'esidu sans "valeur principale"; mais alors, le lemme pr\'ec\'edent permet de conclure.

On suppose maintenant que les hypersurfaces $Y_1,\dots,Y_{p+1}$ sont positives, dans le sens que 
$H^j({\cal E}(\star Y_i))=0$, pour $j>0$, et ${\cal E}$ un faisceau localement libre.
On a alors le th\'eor\`eme suivant:

\begin{thm} \label{th1}
Pour $i\le p$,  le faisceau ${\cal C}^{q,i}_{Y_1\cap \dots\cap Y_i}(\star Y_{i+1})$ est acylique sur $X$. En particulier,
la cohomologie de Dolbeault est calcul\'ee par la suite exacte du lemme pr\'ec\'edent, i.e.:
$$H^i(\Omega^q)\simeq H^0(   {\cal C}^{q,i}_{Y_1\cap \dots\cap Y_i})/\overline\partial(H^0( {\cal C}^{q,i-1}_{Y_1\cap\dots\cap Y_{i-1}}(\star Y_i)))$$
\end{thm}

On a un th\'eor\`eme analogue pour les sous-faisceaux ${\tilde{\cal C}}^{p,q}_Y\subset {\cal C}^{p,q}_Y$ des courants localement r\'esiduels {\it $\partial-$ferm\'es} \`a support dans $Y$.
Si ${\tilde \Omega}^q\subset {\Omega^q}$ est le sous-faisceau des formes holomorphes $\partial-$ferm\'ees, on obtient un isomorphisme analogue:
$$ H^i({\tilde{\Omega}}^q)\simeq H^0(   {\tilde{\cal C}}^{q,i}_{Y_1\cap \dots\cap Y_i})/d(H^0(\overline\partial {\tilde{\cal C}}^{q,i-1}_{Y_1\cap\dots\cap Y_{i-1}}))$$

D'autre part, on peut montrer:
\begin{thm} \label{th2}
1. Tout courant localement localement r\'esiduel $T\in H^0({\cal C}^{q,i-1}_{Y_1\cal\dots\cap Y_{i-1}}(\star Y_{i}))$ (resp. $\partial -$ferm\'e  $T\in H^0({\tilde{\cal C}}^{q,i-1}_{Y_1\cal\dots\cap Y_{i-1}}(\star Y_{i}))$) peut s'\'ecrire comme r\'esidu global:
$$T=Res_{Y_1,\dots,Y_{i-1}}(\Psi)$$, avec $\Psi$ une $q-$forme m\'eromorphe (resp. $\partial-$ferm\'ee) \`a p\^ole dans $Y_1\cap\dots\cap Y_{i}$. 

2. D'autre part, on peut supposer $\Psi$ \`a p\^ole dans $Y_1\cup \dots\cup Y_k$ ssi $T$ est $\overline\partial-$exact (ou, ce qui est \'equivalent, $d-$exact si $T$ est $d-$ferm\'e).
\end{thm}

La deuxi\`eme partie du th\'eor\`eme est d\'emontr\'ee dans \cite{HDS} par des m\'ethodes similaires.

Supposons $X$ une vari\'et\'e projective, et les hypersurfaces $Y_i$ d\'efinies par des diviseurs de Cartier. Alors, la condition de positivit\'e ci-dessus est v\'erifi\'ee, si les fibr\'es associ\'es aux $Y_i$ sont amples, d'apr\`es le th\'eor\`eme d'annulation de Serre.
On en d\'eduit comme corollaire:

\begin{cor}
\label{cor1}
Toute classe de cohomologie $(i,i)$ admet comme repr\'esentant un courant localement r\'esiduel $d-$ferm\'e \`a support dans $Y_1\cap \dots \cap Y_i$.
\end{cor}

D'autre part, supposons $q=n$. On peut encore d\'efinir un sous-faisceau naturel ${{\cal C}^{n,p}}''_{Y}(Y')\subset {\tilde {\cal C}}^{n,p}_Y(\star Y')$, qui sont les sous-faisceaux d\'efinis par les courants de la forme $\omega\wedge[Y]$, avec $\omega$ \`a p\^ole simple sur $Y'$.

On a alors les analogues suivants des deux th\'eor\`emes pr\'ec\'edent.
Le premier th\'eor\`eme est une une variante du th\'eor\`eme principal de \cite{Rosly}.

\begin{thm} \label{th3}
Le sous-complexe 
$$0\to \Omega^n\to {{\cal C}''}^{n,0}( Y_1)\to \dots \to {{\cal C}''}^{n,p}_{Y_1\cap\dots\cap Y_p}\to 0$$
est \'egalement acyclique, de sorte qu'on a un isomorphisme:
$$H^i(\Omega^n)\simeq H^0(   {{\cal C}''}^{n,i}_{Y_1\cap \dots\cap Y_i})/\overline\partial(H^0( {{\cal C}''}^{n,i-1}_{Y_1\cap\dots\cap Y_{i-1}}))$$
\end{thm}

De m\^eme:

\begin{thm} \label{th4}
Pour $i\le p$, tout courant localement r\'esiduel $T\in H^0({{\cal C}''}^{q,i-1}_{Y_1\cap \dots\cap Y_{i-1}}(Y_{i}))$ s'\'ecrit comme r\'esidu global $T=Res_{Y_1,\dots,Y_i}(\Psi)$ d'une $n-$forme m\'eromorphe $\Psi$, \`a p\^oles simples dans $Y_1\cup\dots\cup Y_{i}$. 
De  plus, $T$ est $d$-exact ssi on peut supposer $\Psi$ \`a p\^ole simple dans $Y_1\cup\dots\cup Y_{i-1}$
\end{thm}

On obtient comme corollaire le th\'eor\`eme suivant de \cite{Griffiths}:

\begin{cor}
\label{cor2}
 Supposons $Y_1,\dots,Y_n$ $n$ hypersurfaces positives sur une vari\'et\'e projective $X$, se coupant transversalement en $s$ points distincts $P_1,\dots,P_s$. Soit $c_1,\dots, c_s$ $s$ nombres complexes; une condition n\'ecessaire et suffisante pour qu'il existe une $n-$forme m\'eromorphe $\Psi$, \`a p\^ole simple dans $Y_1\cup \dots \cup Y_n$ admettant comme r\'esidus ponctuels $c_i$ aux points $P_i$, est que $\sum_{i=1}^s{c_i}=0$.
\end{cor}

\section{D\'emonstrations}

\begin{lem}
Le complexe suivant:
$$
\begin{array}{c}
0\to \Omega^q\to \oplus_{1\le j_1\le i}{\Omega(\star Y_{j_1})}\to \oplus_{1\le j_1<j_2\le k}{\Omega(\star Y_{j_1}\cup Y_{j_2})} \\
\to \dots\to \Omega^q(Y_1\cup\dots\cup Y_i)\to {\cal C}^{q,i}_{Y_1\cap\dots\cap Y_i}\to 0
\end{array}$$
o\`u le dernier op\'erateur est donn\'e par le r\'esidu $Res_{Y_1,\dots,Y_i}'$, et les autres op\'erateurs sont les somme altern\'ees du complexe de Cech, est exact
\end{lem}

\dem
En dehors du dernier morphisme, l'exactitude provient de l'exactitude connue du complexe de Cech. La seule particularit\'e ici est que l'on se limite \`a des p\^oles d'ordre fini; mais l'exactitude reste v\'erifi\'ee, car on peut consid\'erer sur un ouvert de Stein $U$ le sous-faisceau des ${\Omega^q}'\subset \Omega^q$ d\'efini pour un ouvert $V\subset U$ par:

${\Omega^q}'(V)$ sont les sections {\it m\'eromorphes} globales, r\'eguli\`ere sur $V$.
D'autre part, l'exactitude au niveau du dernier op\'erateur est un th\'eor\`eme classique de la th\'eorie des r\'esidus (\cite{HDS}).
\CQFD

On peut "tordre" la suite exacte du lemme pr\'ec\'edent en rajoutant des p\^oles sur l'hypersurface $Y^{i+1}$:
$$
\begin{array}{c}
0\to \Omega^q(\star Y_{i+1})\to \oplus_{1\le j_1\le k}{\Omega(\star Y_{j_1}\cup Y_i)}\to \oplus_{1\le j_1<j_2\le i}{\Omega(\star Y_{j_1}\cup Y_{j_2}\cup Y_i)}\\
 \to \dots\to \Omega^q(Y_1\cup\dots\cup Y_{i+1})\to {\cal C}^{q,i}_{Y_1\cap\dots\cap Y_i}(\star Y_{i+1})\to 0
 \end{array}$$

\begin{lem}
Les faisceaux $\Omega^q(\star Y_{j_1}\cup\dots\cup Y_{j_k})$ sont acycliques.
D'apr\`es la suite exacte pr\'ec\'edente, le faisceau
${\cal C}^{q,i}_{Y_1\cap\dots\cap Y_i}(\star Y_{i+1})$ est \'egalemet acyclique.
\end{lem}
\dem
On sait d'apr\`es l'hypoth\`ese de positivit\'e que si $L_j$ est le fibr\'e associ\'e \`a $Y^i$, on a:
$H^j(\Omega\otimes L_{j_1}^{s_1}\otimes L_{j_k}^{s_k})=0$ pour $j>0$, et $s_1,\times s_k$ assez grands. Par passage \`a la limite, on en d\'eduit que $H^j(\Omega(\star (Y_{j_1}\cup\dots\cup Y_{j_k}))=0$ pour $j>0$.
D'autre part, l'acyclicit\'e de ${\cal C}^{q,i}_{Y_1\cap\dots\cap Y_i}(\star Y_{i+1})$ d\'ecoule d'apr\`es la suite exacte pr\'ec\'edente qu'un faisceau ayant une "r\'esolution" acyclique est lui-m\^eme acyclique.
\CQFD

En particulier, comme le foncteur $\Gamma$ des sections globales est exact dans la cat\'egorie des faisceaux acycliques, on en d\'eduit, en appliquant $\Gamma$ \`a la suite exacte pr\'ec\'edente, que 

{\it Le morphisme $Res_{Y_1,\dots,Y_{i}}(\Omega^q(\star Y_1\cup\dots\cup Y_{i+1})\to H^0({\cal C}^{q,i}_{Y_1\cap \dots\cap Y_i}(\star Y_{i+1})$ est surjectif,}
ce qui d\'emontre la premi\`ere partie du th\'eor\`eme \ref{th2}.

Soit $T$ le courant $T:=Res_{Y_1,\dots,Y_i}(\Psi)$, avec $\Psi$ m\'eromorphe \`a p\^oles dans $Y_1\cup\dots\cup Y_{i+1}$. Si $\Psi$ est \`a p\^ole contenu dans $Y_1\cup\dots\cup Y_{i}$, on sait que $T=\overline\partial Res_{Y_1,\dots,Y_{i-1}}(\Psi)$, et dont $T$ est $\overline\partial-$exact.
D'autre part, si $T$ est $\overline\partial-$exact, on en d\'eduit d'apr\`es le th\'eor\`eme \ref{th1} que $T$ s'\'ecrit $\overline\partial T'$, avec $T'$ un courant de $H^0({\cal C}^{q,i-1}_{Y_1\cap \dots\cap Y_{i-1}}(\star Y_i))$; d'apr\`es ce qu'on vient de montrer, pour le degr\'e $i-1$, on voit que $T'$ s'\'ecrit $T'=Res_{Y_1,\dots,Y_{i-1}}(\Psi')$, avec $\Psi'$ \`a p\^ole contenu dans $Y_1\cup\dots\cup Y_{i-1}$; on a donc $T=\overline\partial T'=Res_{Y_1,\dots,Y_i}(\Psi')$, ce qui termine la d\'emonstration de la deuxi\`eme partie du th\'eor\`eme \ref{th2}.

Si $X$ est une vari\'et\'e projective, et $Y_i$ l'annulation d'une section globale d'un fibr\'e lin\'eaire ample, la condition de positivit\'e sur les $Y_i$ est v\'erifi\'ee d'apr\`es le th\'eor\`eme d'annulation de Serre. Le corollaire \ref{cor1} s'en d\'eduit imm\'ediatement.

D'autre part, la d\'emonstration du th\'eor\`eme \ref{th3} se fait de la mani\`ere que celle du th\'eor\`eme \ref{th2}, l'acyclicit\'e des faisceaux ${\Omega}^n(Y_{j_1}\cup Y_{j_k})$ se d\'emontrant cette fois par le th\'eor\`eme d'annulation de Kodaira.
Pour le corollaire \ref{cor2}, l'annulation $\sum_{i=1}^s{c_s}=0$ implique que le courant d'"\'evaluation" $\sum_{i=1}^s{c_i}[P_i]$ est $d-$exact, et s'\'ecrit donc comme un r\'esidu global d'apr\`es le th\'eor\`eme \ref{th4}.

D'apr\`es ce qui pr\'ec\`ede, une nouvelle formulation de la conjecture de Hodge devient:

{\it Si un courant localement r\'esiduel $T$ de bidegr\'e $(i,i)$, $d-$ferm\'e et \`a support dans $Y_1\cap \dots\cap Y_i$, est tel que pour toute sous-vari\'et\'e complexe de dimension $i$ $C$, la somme des r\'esidus $T.C$ est nulle, de sorte que le courant r\'esiduel correspond sur $C$ est globalement r\'esiduel, alors si la classe de coholomogie de $T$ est enti\`ere, on peut en d\'eduire que $T$ est globalement r\'esiduel.}

Cette nouvelle formulation serait int\'eressante \`a \'etudier sur des exemples particuliers.

{{\bf Remerciements.}
Je remerciement vivement Sebastien Boucksom pour sa collaboration, par de multiples discussions instructives et des critiques constructives. En particulier, pour m'avoir signal\'e l'importance du premier th\'eor\`eme \ref{th1}, et la r\'ef\'erence \cite{Rosly}, qui avec \cite{HDS} est une des principales sources de ce travail.
}

\end{document}